\documentclass[11pt]{pjm} 

\usepackage{epsfig}
  
\title[Product-to-sum formula]{A
product-to-sum formula for the quantum group of $SL(2,{\mathbb C})$} 
  
\author{R{\u{a}}zvan Gelca} 
\address{Department of Mathematics and Statistics, 
Texas Tech University, Lubbock, TX 79409 and Institute of Mathematics
of the Romanian Academy, Bucharest, Romania}
\email{rgelca@math.ttu.edu}
  
\newtheorem{thm}{Theorem}[section] 
 
\newtheorem{lem}[thm]{Lemma}

\theoremstyle{definition}

\theoremstyle{remark} 
  
\newtheorem{rem}[thm]{remark}

\begin{document} 
  
\begin{abstract} 
  
This paper  exhibits a  product-to-sum formula for the observables of a 
certain  quantization 
of the moduli space of flat $SU(2)$-connections on the torus.
This quantization was  defined using the topological quantum field theory 
that was developed by  Reshetikhin and Turaev   from 
the quantum group of $SL(2,{\mathbb C})$
at roots of unity. As a corollary it is shown  that the 
algebra of quantum observables is a subalgebra of the noncommutative
torus with  rational rotation angle. 
The proof  uses topological quantum field theory with
corners, and is based on the  description of   the matrices of the observables
in a canonical basis of the Hilbert space of the quantization.
  
\end{abstract} 
  
\maketitle 
  
\section{Introduction}\label{intro} 
  
 In this paper we  prove  a product-to-sum formula that holds for the 
quantization of the moduli space of flat $SU(2)$-connections on the torus.
A quantization of this moduli space was given by Witten in \cite{Wi}
in conjunction with the Jones polynomial, using a path integral
on the space of connections. Witten's 
path integrals were put into rigorous framework 
 by Reshetikhin and Turaev in \cite{RT} using
quantum groups. This paper is based on  their construction.

The algebra of observables is densely generated by the functions
of the form $\cos 2\pi (px+qy)$, $p,q$ integers.
 The quantization of such a function
consists of the coloring of the curve of slope $p'/q'$ on the torus
by the virtual representation $V^{n+1}-V^{n-1}$, where 
$n$ is the greatest common divisor of $p$ and $q$, $p'=p/n$, $q'=q/n$,
and we denote by $V^{k}$ the $k$-dimensional irreducible 
representation of the quantum group of $SL(2,{\mathbb C})$. 
In the work Reshetikhin and Turaev  $V^k$ is only defined for
$k=1,2,\ldots , r-1$, but we will see below how this  definition
is  extended to all integers. Let us denote by $C(p,q)$ the
operator that is associated by the quantization to
$\cos 2\pi (px+qy)$.

A product-to-sum formula for the version of these operators  defined
using  the Kauffman bracket was discovered in \cite{FG} and that is what 
inspired the following result.
  
\begin{thm}\label{mainresult}
In any level $r$ and for  any integers $m,n,p,q$ the following 
product-to-sum formula holds
\begin{eqnarray*}
C(m,n)*C(p,q)=t^{|^{mn}_{pq}|}C(m+p,n+q)+t^{-|^{mn}_{pq}|}C(m-p,n-q).
\end{eqnarray*}
\end{thm}

So we see that the quantization of the algebra of
the observables on the moduli space of flat $SU(2)$-connections on
the torus is isomorphic to a subalgebra of the noncommutative torus
defined by Rieffel \cite{rieffel}.
To be more rigorous, the above product-to-sum formula
induces a $*$-product on the algebra of functions on the character
variety on the torus. This $*$-algebra  is a subalgebra
of the noncommutative torus.

  The noncommutative torus is
a $C^*$-algebra generated by $U$, $V$, $U^{-1}$, $V^{-1}$ where $U$
and $V$ are two unitary operators satisfying the exponential form
of the Heisenberg noncommutation relations $UV=t^2VU$.
The inclusion  is given by
\begin{eqnarray*}
C(p,q)\rightarrow t^{-pq}(U^pV^q+U^{-p}V^{-q}).
\end{eqnarray*}  
The elements $\frac{1}{2}t^{-pq}(U^pV^q+U^{-p}V^{-q})$ are
the noncommutative cosines, so  Theorem \ref{mainresult} is 
nothing but the product-to-sum formula for noncommutative cosines. 
Of course, $t$  root of unity is  the situation 
Rieffel did not consider, for here  the rotation angle
in the definition of the noncommutative torus is rational. 

What we consider to be a nice feature of this paper is that the proof
of this  result about the quantization of moduli spaces 
uses topological cut-and-paste techniques. 

 \section{The background of the problem}\label{background} 

Let us briefly discuss the case of the torus and then state the main result.
The moduli space of flat $SU(2)$-connections on the torus is
the ``pillow case'', the quotient of the complex plane by the lattice
${\mathbb Z}[i]$ and by the symmetry with respect to the origin. 
This space is the same as the character variety of $SU(2)$-representations
of the fundamental group of the torus. 

The algebra of observables, i.e. the algebra  of functions on this
space, is generated for example
by the functions $\sin 2\pi n(px+qy)/\sin 2\pi (px+qy)$, where $n,p,q$ are
integer numbers with $n\geq 0$ and $p,q$ coprime. 
A function of this form has  a geometric interpretation, namely it associates
to a conjugacy class of flat connections   
the trace of the holonomy of the connection along the curve of
slope $p/q$, where the trace is  computed in the  $n$-dimensional irreducible
representation of $SU(2)$.    
Witten defined the path integral
 only for functions of this form, but the 
definition can be extended to any (smooth) function on the character variety
to yield an operator 
\begin{eqnarray*}
op(f)=\int e^{i{\mathcal L}(A)}f(A){\mathcal D}A
\end{eqnarray*}
where $f(A)$ is defined using  approximations of  $f$ by sums of
traces of holonomies of the type describe above,  
\begin{eqnarray*}
{\mathcal L}(A)=
\frac{k}{4\pi}Tr\int_{{\mathbb T}^2\times I}\left(A\wedge dA +\frac{2}{3}A
\wedge A\wedge A\right)
\end{eqnarray*}
is the
 Chern-Simmons functional for connections on the cylinder over the torus, 
$k$ being the level of quantization, 
and the path integral is taken over all connections $A$ that interpolate
between two flat connections on the  boundary components  of the cylinder
over the torus.

Witten's construction was made rigorous by Reshetikhin and Turaev using
the quantum group of $SL(2,{\mathbb C})$ at roots of unity.
Here is how this is done for the torus.
Fix a level $r\geq 3$, and let $t=e^{\frac{\pi i }{2r}}$.
 For an integer $n$ define 
$[n]=(t^{2n}-t^{-2n})/(t^2-t^{-2})$. The quantum analogue of the universal
enveloping algebra of $sl(2,{\mathbb C})$ is the algebra
${\mathbb U}_t$ with generators 
$X,Y,K, \bar{K}$ satisfying
\begin{eqnarray*}
& & \bar{K}=K^{-1}, \quad  KX=t^2XK, \quad  KY=t^{-2}YK\\
& & XY-YX=\frac{K^2-K^{-2}}{t^2-t^{-2}}, \quad
 X^r=Y^r=0, \quad K^{4r}=1.
\end{eqnarray*}

This is a Hopf algebra, thus its representations form a ring
under direct sum and tensor product. Modulo a small technicality
where one  factors by the  part of quantum trace zero,
this ring contains a subring generated by finitely many
irreducible representations $V^1,V^2, \ldots , V^{r-1}$.
The representation $V^k$ has the  basis
$e_{-\frac{k-1}{2}},
e_{-\frac{k-3}{2}}\ldots , e_{\frac{k-1}{2}}$ and ${\mathbb U}_t$ acts
on this basis by
\begin{eqnarray*}
& & Xe_j=[m+j+1]e_{j+1}\\
 & & Ye_j=[m-j+1]e_{j-1}\\
& & Ke_j=t^{2j}e_j.
\end{eqnarray*}
Using the quantum version of the Clebsch-Gordon theorem, 
these irreducible representations can be defined recursively
by $V^{n+1}=V^2\otimes V^{n}-V^{n-1}$. This definition works
only for $2\leq n\leq r-2$.

With these  representations we color any   knot $K$ 
in a 3-dimensional
manifold, and the  coloring is denoted by $V^n(K)$.
Let $\alpha $ be the  core $S^1\times \{0\}$ of 
 the solid torus $S^1\times {\mathbb D}$, ${\mathbb D}=
\{ z, |z|\leq 1\}$. 
The Hilbert space of the torus $V({\mathbb T}^2)$
has an orthonormal basis given
by $V^{k}(\alpha)$, $k=1,2,\ldots, r-1$.  Now the convention is made
that $V^{r} =0$ and with this in mind we can extend the 
definition of $V^n(\alpha )$ to all $n$ via the recursion
\begin{eqnarray*}
V^{n+1}(\alpha )=V^2(\alpha )\otimes V^{n}(\alpha )-V^{n-1}(\alpha).
\end{eqnarray*}
Here $V^2(\alpha )\otimes V^{n}(\alpha)$ is just $\alpha $ colored by
the representation $V^2\otimes V^n$. 

To make sense of these notations  recall the   pairing $<\cdot , \cdot >$ 
on $V({\mathbb T}^2)$ obtained by gluing two solid tori such that
the meridian of the first is identified with the longitude
of the second and vice versa, as to obtain a 
 3-sphere. This pairing is not an inner product but is nondegenerate.
Pairing two elements of $V({\mathbb T}^2)$ yields a colored link
in the 3-sphere, and Reshetikhin and Turaev described in \cite{RT} a
way of associating a number to this link using the ribbon algebra structure
of ${\mathbb  U}_t$ (see also \cite{KR}, \cite{Tu}).
 This number is the value of the pairing. 
For example
it is well known that 
\begin{eqnarray*}
<V^k(\alpha),V^m(\alpha)>=[mn].
\end{eqnarray*}

Let $p$ and $q$ be two integers, $n$ their greatest common divisor
and $p'=p/n$, $q'=q/n$.
We describe now the operator $S(p,q)$ that corresponds to
the function $\sin 2\pi n(p'x+q'y)/\sin 2\pi (p'x+q'y)$ ($S$ stands for 
sine). Since the pairing is nondegenerate, it suffices to 
describe $<S(p,q)V^k(\alpha), V^m(\alpha)>$. This number is computed 
 by placing the curve of slope $p'/q'$ on the
boundary of the solid torus, then gluing another torus to it
to obtain the 3-sphere, coloring the $p'/q'$-curve by $V^n$, the
core of the first torus by $V^k$ and the core of the second torus
by $V^m$ and then evaluating the colored link diagram as in 
\cite{RT}. 
We neglected  the discussion about the 
framing, and the invariants of Reshetikhin and Turaev are defined for framed
knots and links. The curve $\alpha$ will always be  framed  by the
annulus $S^1\times [0,1]$ while the curve of slope $p'/q'$ has
the vertical framing (i.e. the framing defined by the vector field 
orthogonal to the torus). 

The quantization of the function $f(x,y)=2\cos 2\pi(px+qy)$ is the operator
\begin{eqnarray*}
C(p,q)=C(np',nq')=S((n+1)p', (n+1)q')-S((n-1)p',(n-1)q').
\end{eqnarray*}
Here $C$ stands for cosine and the definition 
is motivated by the trigonometric formula
\begin{eqnarray*}
2\cos nx=\frac{\sin (n+1)x}{\sin x}-\frac{\sin(n-1)x}{\sin x}.
\end{eqnarray*}
These are the operators for which we prove the product-to-sum formula.
  
\section{Proof of the main result}\label{proofofresult} 

  Unlike the case of \cite{FG}, here an easy inductive proof will not
work. The proof we give  below  was inspired by the computation of
the colored Jones polynomials of torus knots  from
\cite{gelca2}. It uses  an apparatus called {\em topological quantum field
theory with corners}, which  enables the computation of quantum invariants
of 3-manifolds through a successive application of
 axioms. This system of axioms is
in spirit analogous to the   Eilenberg-MacLane system of axioms for homology.
The fundamental principles of a TQFT with corners were summarized  by Walker
in \cite{walker}. The construction of the basic data
for the case of the quantum group 
${\mathbb U}_t$ was initiated in \cite{FK} and completed in \cite{gelca1}.

Let us discuss just the  facts about the $SL(2,{\mathbb C})$ TQFT with corners
 needed in the proof of Theorem \ref{mainresult}. 
As explained in \cite{atiyah}, a TQFT in dimension three consists of
a modular functor $V$ that associates vector spaces to 
surfaces and isomorphisms to homeomorphisms, and a partition function
that associates to each three dimensional manifold a vector in the 
vector space of its boundary. 

 Walker refined this point of view
in \cite{walker}. He considered decompositions of surfaces into disks,
annuli, and pairs of pants, called DAP-decompositions, which 
correspond through the functor to  basis' of the vector space.
His DAP-decompositions (which should probably be called rigid
structures) involve more structure than that. This structure consists of 
some curves, called seams and drawn like dotted lines in diagrams,
which keep track of the twistings that occur and are not detected   
by the change in DAP-decompositions. Transformations of DAP-decompositions
are called moves and have the vector space correspondent of 
the change of basis. They are useful in preparing the boundary
of a 3-manifold for a gluing with corners. Let us point out that
for ${\mathbb U}_t$ one needs also an
orientation of the decomposition curves, as explained  in
\cite{gelca1} where it was  encoded using the Klein four group, but this
fact is irrelevant for  all computations below, so we do not
refer to it again. 

Also, Walker considered framed, or extended manifolds, which are
 pairs $(M,n)$ with $M$ a manifold and $n$ an integer,
necessary to make the partition function well defined.
In this  paper we begin the computation with manifolds with
framing zero. The framing only changes when we perform 
moves on boundary tori, and there the change is computed using  the
Shale-Weil cocycle. However all relevant changes are performed back and
forth, so the framing cancels. Thus we simply ignore it. 

Following \cite{walker} we denote $X=\sqrt{\sum_{j=1}^{r-1}[j]^2}$
(not to be confused with the generator of ${\mathbb U}_t$).
The vector space of an annulus has basis $\beta_j^j$, $j=1,2,\ldots, r-1$,
with  pairing given by 
\begin{eqnarray*}
<\beta_j^j,\beta_k^k>=\delta_{j,k}
X/[j].
\end{eqnarray*}   
Gluing the boundaries of an annulus we obtain a torus, with the same
basis for the vector space.
The vector space of the torus $V({\mathbb T}^2)$ 
 is the same as the 
Hilbert space of the quantization of the moduli space of flat 
$SU(2)$-connections on the torus. 
 The pairing on the torus makes $\beta_j^j$ an
orthonormal basis, but we don't need this fact now. 
The moves  $S$ and $T$ on the torus are described in Figure 1.
The 
 $(m,n)$-entry of the matrix of $S$ is $[mn]$. The move $T$ is diagonal
and its $j$th entry is $t^{j^2-1}$. 

\begin{figure}[htbp]
\centering
\leavevmode
\epsfxsize=4.4in
\epsfysize=1.1in
\epsfbox{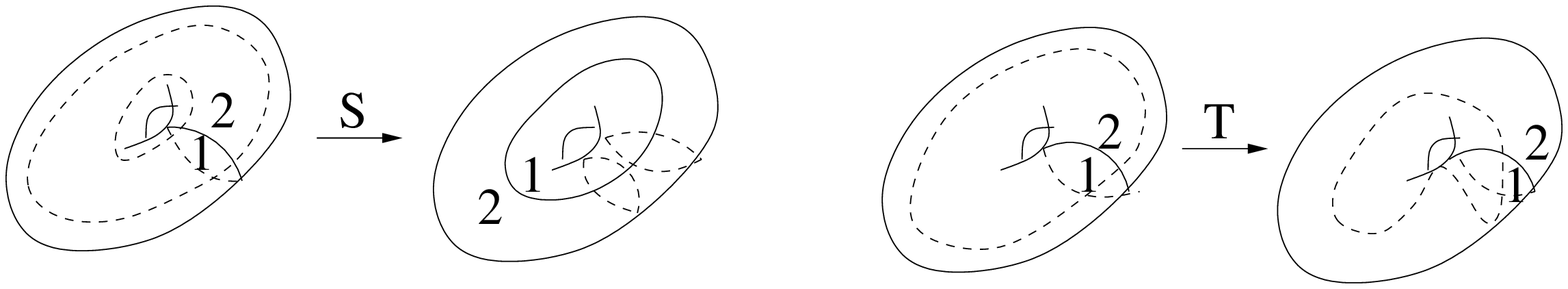}

Figure 1.    
\end{figure}

The quantum invariant (i.e. partition function)
of the cylinder over a surface is the identity
matrix. So the quantum invariant of 
the cylinder over an annulus is $\sum _{n=1}^{r-1}\frac{[n]}{X}\beta_n^n\otimes
\beta_n^n\otimes \beta_n^n\otimes\beta_n^n$ (when taking the cylinder
over the annulus the boundary of the solid torus will be canonically
decomposed into four annuli).
If the DAP-decomposition of a  3-manifold involves two disjoint annuli,
then its  invariant can be written in the form
\begin{eqnarray*}
\sum _{k,j=1}^{r-1} \beta_k^k\otimes\beta_j^j \otimes v_{k,j}.
\end{eqnarray*} 
 Gluing the two annuli produces a 3-manifold that in the newly
obtained DAP-decomposition has the invariant equal to
$\sum _{k=1}^{r-1}\frac{X}{[k]}v_{k,k}$.

For the proof of the theorem we need  the following formula, 
which can be  checked using the 
fact that the sum of the roots of unity is zero. 

\begin{lem}\label{rootsofunity}
Let $a$, $b$, $c$, $d$ and $e$ be integers. Then
\begin{eqnarray*}
\sum_{x,y=1}^{r-1} [ax]t^{bx^2}[cy]([x(y+d)]t^{2ey}+[x(y-d)]t^{-2ey})=\\
X^2t^{bc^2+be^2-2de}([a(c+e)]t^{2(be-d)c}+[a(c-e)]t^{-2(be-d)c}).
\end{eqnarray*}
\end{lem}

The proof of the main result is based on a theorem  that is of
interest in itself. 
It describes how the observable $C(p,q)$ acts on the 
Hilbert space of the quantization. With  our convention
for defining  
$V^k({\alpha})$ for all $k$ we have  

\begin{thm}
In any level $r$ and for any integers $p,q$ and $k$ the following formula
holds
\begin{eqnarray*}
C(p,q)V^k(\alpha)=t^{-pq}\left(t^{2qk}V^{k-p}(\alpha)+t^{-2qk}V^{k+p}(\alpha)
\right).
\end{eqnarray*}
\end{thm}

\begin{proof}
It suffices to check that
the two sides of the equality yield the same results when paired with
all $V^m(\alpha)$. We must show  that  
\begin{eqnarray*}
& & <C(p,q)V^k(\alpha),V^m(\alpha)> 
=\frac{t^{-pq}}{t^2-t^{-2}}\\ & & \quad \times \left(
t^{2(qk-pm+km)}-t^{2(qk+pm-km)}+t^{2(-qk+pm+km)}-t^{2(-qk-pm -km)}\right)
\end{eqnarray*}

Let $d$ to be the greatest common divisor of 
$p$ and $q$, $p'=p/d$ and $q'=q/d$
We concentrate  first on the computation of 
\begin{eqnarray}\label{invariant1}
<S((d+1)p',(d+1)q')V^k(\alpha),V^{m}(\alpha)>
\end{eqnarray}
and 
\begin{eqnarray}\label{invariant2}
<S((d-1)p',(d-1)q')V^k(\alpha),V^{m}(\alpha)>
\end{eqnarray}
Here is the place where we use the topological quantum field
theory with corners from  \cite{gelca1}.

The  expression in (\ref{invariant1}) is the invariant
of the link that has one component equal to the $(p',q')$-curve on a torus 
colored by $V^{d+1}$, and the other two components the cores
of the two solid tori that lie on one side and  the other
of the torus knot, colored by $V^k$ (the one inside) and
$V^m$ (the one outside). The expression in (\ref{invariant2}) is the invariant
of the same link but with the $(p',q')$-curve colored by $V^{d-1}$. 
It was shown in \cite{gelca2} that this number is equal to $X^{-1}$
times 
the coordinate of $\beta _{d+1}^{d+1}d\otimes \beta_k^k\otimes \beta _m^m$
of the vector that is the quantum invariant of the link complement.



Let us produce the complement of this link by gluing together two
simple 3-manifolds, whose quantum invariants are easy to compute.
Consider first the cylinder over an annulus $A$ and glue its
ends to obtain the manifold $A\times S^1$. In the basis of the vector
space of $V({\mathbb T}^2\times {\mathbb T}^2)$  determined 
by the DAP-decomposition
$\partial A\times \{1\}$ the invariant of this manifold is
$\sum _{k}  \beta _k^k\otimes \beta_k^k$.
 Take another copy of the same manifold.
Change the decomposition curves  of the exterior torus  of the first 
manifold to the $p'/q'$-curve and of the exterior torus of the second 
manifold to the longitude. 
Of course, to do this on  the second manifold  we  
apply the $S$-move and so the invariant of the second manifold 
changes to
\begin{eqnarray*}
\frac{1}{X}\sum _{\delta,j_{n+1}}[dj_{n+1}]\beta_\delta^\delta \otimes 
\beta _{j_{n+1}}^{j_{n+1}}.
\end{eqnarray*}
With the first manifold the story is more complicated. Consider the continued
fraction expansion
\begin{eqnarray*}
\frac{q'}{p'}=-\frac{1}{
-a_1-\frac{1}{
-a_2-\cdots \frac{1}{-a_n}}}.
\end{eqnarray*}
The required move on the boundary is then
$ST^{-a_n}ST^{-a_{n-1}}S\cdots ST^{-a_1}S$.
So the invariant of the first manifold in the new
DAP-decomposition is
\begin{eqnarray*}
X^{-n-1}\sum_{j_1,\cdots ,j_{n+1}} [j_{n+1}j_n]t^{-a_n(j_n^2-1)}[j_nj_{n-1}]
\cdots [j_2j_1]t^{-a_1(j_1^2-1)}[j_1k]\beta_k^k\otimes 
\beta _{j_{n+1}}^{j_{n+1}}.
\end{eqnarray*}

Now expand one annulus in the exterior tori  of each of the two manifolds.
Then glue just one annulus from the the first manifold to one
annulus from the second. This way we obtain the complement of the
link in discussion. One of its boundary tori is decomposed into two
annuli. Contract one of them. Since the gluing introduces
a factor of $X/[j_{n+1}]$, the invariant of the manifold is 
\begin{eqnarray*}
X^{-n-1}\sum_{j_1,\cdots ,j_{n+1},\delta}\frac{[\delta j_{n+1}]}{[j_{n+1}]}
 [j_{n+1}j_n]t^{-a_n(j_n^2-1)}[j_nj_{n-1}]
\cdots [j_2j_1]t^{-a_1(j_1^2-1)}\\
\times [j_1k]\beta _{j_{n+1}}^{j_{n+1}}\otimes 
\beta_\delta^\delta\otimes \beta_k^k.
\end{eqnarray*}
At this moment we have the right 3-manifold but with the wrong 
DAP-decomposition. We need to fix the DAP-decomposition of the torus
that corresponds to the basis element $ \beta _{j_{n+1}}^{j_{n+1}}$ 
(the boundary of the regular neighborhood of the $(p',q')$-curve)
such as to transform the decomposition curve into the meridian of the link
component. For this we apply the move  
$(ST^{-a_n}ST^{-a_{n_1}}S\cdots ST^{-a_1}S)^{-1}$. We obtain the following
expression for the invariant of the extended manifold 
\begin{eqnarray*}
X^{-2n-2}\sum_{j_1,\ldots ,j_{2n+2},\delta,k,m}[mj_{2n+2}][j_{2n+2}j_{2n+1}]
t^{a_1j_{2n+1}^2}\cdots [j_{n+2}j_{n+1}]\frac{[\delta j_{n+1}]}{[j_{n+1}]}
\\ \times  [j_{n+1}j_n]
 t^{-a_nj_n^2}[j_nj_{n-1}]
\cdots 
[j_2j_1]t^{-a_1j_1^2}[j_1k]\beta_\delta^\delta 
\otimes\beta_k^k\otimes \beta_m^m
\end{eqnarray*}
(in this formula we already  reduced $t^{a_k}$ and $t^{-a_k} $, $1\leq 
k\leq n$). 

It is important to observe that after performing the 
described operations the seams came right, so no further twistings
are necessary.

Now fix $k$ and $m$, let $\delta =d\pm 1$ and focus on the coefficients of  
$\beta_{d\pm 1}^{d\pm 1}\otimes\beta_k^k\otimes \beta_m^m$. 
Multiplied by $X$ these are 
the colored Jones polynomial of the link \cite{gelca2} with the 
$(p',q')$-curve colored by the $d+1$-, respectively $d-1$-dimensional 
irreducible representation of ${\mathbb U}_t$. 
Since $C(p,q)=S((d+1)p',(d+1)dq')-S((d-1)p',(d-1)q')$
and also 
\begin{eqnarray*}
\frac{[(d+1)j_{n+1}]}{[j_{n+1}]}-\frac{[(d-1)j_{n+1}]}{[j_{n+1}]}=
t^{2j_{n+1}}+t^{-2j_{n+1}},
\end{eqnarray*}
 we deduce that
 the value of $<C(p,q)V^k(\alpha),V^{m}(\alpha)>$
is equal to 
\begin{eqnarray*}
X^{-2n-1}\sum_{j_1,\ldots ,j_{2n+2}}[mj_{2n+2}][j_{2n+2}j_{2n+1}]
t^{a_1j_{2n+1}^2}\cdots [j_{n+2}j_{n+1}](t^{2dj_{n+1}}+t^{-2dj_{n+1}})
 \\
\times [j_{n+1}j_n]t^{-a_nj_n^2}[j_nj_{n-1}]
\cdots 
[j_2j_1]t^{-a_1j_1^2}[j_1k].
\end{eqnarray*}
We want to compute these iterated Gauss sums.
We apply successively  Lemma \ref{rootsofunity}
 starting with $x=j_n, y=j_{n+1}$, then
$x=j_{n-1}, y=j_{n+2}$ and so on to obtain 

\begin{eqnarray*}
& & X^{-2n+1}\sum_{j_1,\ldots ,j_{2n+2}}t^{-a_nd^2}
[mj_{2n+2}][j_{2n+2}j_{2n+1}]
t^{a_1j_{2n+1}^2}\cdots [j_{n+3}j_{n+2}]\\
& & \times ([j_{n-1}(j_{n+2}+d)]t^{-2a_ndj_{n+2}}
+[j_{n-1}(j_{n+2}-d)]t^{2a_ndj_n+2}
)\\
& & \times t^{-a_{n-1}j_{n-1}^2}[j_{n-1}j_{n-2}]
\cdots [j_2j_1]t^{-a_1j_1^2}[j_1k]
\end{eqnarray*}
\begin{eqnarray*}
& & =X^{-2n+3}
\sum_{j_1,\ldots ,j_{2n+2}}t^{-a_n(a_na_{n-1}-1)d^2}
[mj_{2n+2}][j_{2n+2}j_{2n+1}]
t^{a_1j_{2n+1}^2}\cdots [j_{n+4}j_{n+3}]\\
& & \times
([j_{n-2}(j_{n+3}+a_nd)]t^{-2(a_na_{n-1}-1)dj_n+3}+[j_{n-2}(j_{n+3}-a_nd)]
t^{2(a_na_{n-1}-1)dj_n+3}
)\\
& & \times t^{-a_{n-2}j_{n-2}^2}
[j_{n-2}j_{n-3}]
\cdots 
[j_2j_1]t^{-a_1j_1^2}[j_1k]=\cdots 
\end{eqnarray*}
\begin{eqnarray*}
& & =
X^{-1}\sum_{j_{2n+1},j_{2n+2}}t^{-p'q'd^2}
[mj_{2n+2}][j_{2n+2}j_{2n+1}]\\
& & \times ([kj_{2n+1}+kdq']t^{-2dp'j_{2n+1}}
+[kj_{2n+1}-kdq']
t^{2dp'j_{2n+1}})\\
& & = X^{-1}
t^{-pq}\sum_{j_{2n+1},j_{2n+2}=1}^{r-1}
[mj_{2n+2}][j_{2n+2}j_{2n+1}]\\ & & \times 
([kj_{2n+1}+kq]t^{-2pj_{2n+1}}
+[kj_{2n+1}-kq]
t^{2pj_{2n+1}})
\end{eqnarray*}

This sum is equal to
\begin{eqnarray*}
t^{-pq}([k(m+q)]t^{-2mp}+[k(m-q)]t^{2mp})
\end{eqnarray*}
and the theorem is proved. 
\end{proof}

Let us  see how Theorem \ref{mainresult} follows from this result.
We compute:
\begin{eqnarray*}
& & C(m,n)*C(p,q)V^k(\alpha)\\
& & \quad =C(m,n)t^{-pq}(t^{2kq}V^{k-p}(\alpha)+
t^{-2kq}V^{k+p}(\alpha))\\
& & \quad = t^{-pq-mn}(t^{2kq+2(k-p)n}V^{k-p-m}(\alpha)+t^{2kq-2(k-p)n}V^{k-p+m}
(\alpha)\\
& & \quad +
t^{-2kq+2(k+p)n}V^{k+p-m}(\alpha)+t^{-2kq-2(k+p)n}V^{k+p+m}(\alpha)).
\end{eqnarray*}
Also  
\begin{eqnarray*}
& & C(m+p,n+q)V^k(\alpha)\\
& & \quad =t^{-(m+p)(n+q)}(t^{2k(n+q)}V^{k-m-p}(\alpha)
+t^{-2k(n+q)}V^{k+m+p}(\alpha))
\end{eqnarray*}
on the one hand, and 
\begin{eqnarray*}
& & C(m-p,n-q)V^k(\alpha)\\
& & \quad = t^{-(m-p)(n-q)}(t^{2k(n-q)}V^{k-m+p}(\alpha)
+t^{-2k(n-q)}V^{k+m-p}(\alpha))
\end{eqnarray*}
on the other.
An easy check of coefficients shows that the product-to-sum formula holds
on the basis of the Hilbert space of the quantization, and we are done.

\begin{rem}
The same method can be applied to prove the product-to-sum formula
for the Kauffman bracket skein algebra of the torus (see \cite{FG})
using the topological quantum field theory with corners constructed
in \cite{gelca3}. 
\end{rem}
We conclude with a question. Is the product-to-sum formula true
when the deformation variable $t$ is not a root of unity? 
The quantization for arbitrary $t$ is described in \cite{AS}, and
a proof of the formula should involve the analysis of the 
quasitriangular $R$ matrix of ${\mathbb U}_t$.

\end{document}